\documentclass[11pt]{amsart}
\usepackage{latexsym}
\usepackage{amsfonts}
\usepackage{amsthm}
\usepackage{graphicx}
\usepackage{amssymb}
\usepackage{amsmath}
\usepackage{amssymb}
\usepackage{color}

\author[]{Vladimir Shpilrain}
\address{Department of Mathematics, The City  College  of New York, New York,
NY 10031} \email{shpil@groups.sci.ccny.cuny.edu}
\thanks{Research of the author was partially supported by
the NSF grant DMS-0914778. }


\dedicatory{Dedicated to Paul Schupp in appreciation of his
contributions to mathematics and computer science}

\newtheorem{proposition}{Proposition}[section]

\newtheorem{example}{Example}
\newtheorem{problem}{Problem}

\newtheorem{lemma}{Lemma}
\newtheorem{thm}{Theorem}[section]
\newtheorem{prob}[thm]{Problem}
\newtheorem{defn}[thm]{Definition}

\begin{document}

\title[Sublinear time algorithms]{Sublinear time algorithms in the theory of groups and
semigroups}

\begin{abstract}
Sublinear time algorithms represent a new paradigm in computing,
where an algorithm must give some sort of an answer after inspecting
only a small portion of the input.  The most typical situation where
sublinear time algorithms are considered is {\it property testing}.
There are several interesting contexts where one can test properties
in sublinear time. A canonical example is graph colorability. To
tell that a given graph is not $k$-colorable, it is often sufficient
to inspect just one vertex with incident edges: if the degree of a
vertex is greater than $k$, then the graph is not $k$-colorable.

It is a challenging and interesting task to find algebraic
properties that could be tested in sublinear time. In this paper, we
address several algorithmic  problems in the theory of groups and
semigroups that may admit sublinear time solution, at least for
``most" inputs.
\end{abstract}

\maketitle

\section{Introduction}

 Typically, to give some information about an input, an algorithm
 should at least ``read" the entire input, which takes linear time in
 ``length", or complexity, of the latter. Thus, linear time was
 usually considered the ``golden standard" of achievement in computational
 complexity theory.

Sublinear time algorithms represent a new paradigm in computing,
where an algorithm must give some sort of an answer after inspecting
only a small portion of the input. Given that reading some data
takes too long, it is natural to ask what properties of the data can
be detected by sublinear time algorithms that read only a small
portion of the data. Sublinear time algorithms for {\it decision
problems} are examples of {\it  property testing} algorithms.

 In broad terms, property testing is the study of
the following class of problems:

\begin{quote}
Given the ability to perform local queries concerning a particular
object (e.g., a graph, or a group element), the task is to determine
whether or not the object has a specific property. The task should
be performed by inspecting only a small (possibly randomly selected)
part of the whole object.
\end{quote}

Often, a small probability of failure is allowed, especially when
efficiency is more important than accuracy; this makes a difference
with ``usual" decision algorithms that have to give correct answers
for all inputs. (By ``failure" here we mean a situation where an
algorithm cannot give a conclusive answer, but we do not allow an
algorithm to give a wrong answer.) In this sense, one of the ideas
behind using sublinear time algorithms is similar to that of using
{\it genericity}, i.e., assessing complexity of an algorithm on
``most" inputs, see e. g. \cite{KMSS, KS, KSS}.

Property testing algorithms offer several benefits: they save time,
are good in settings where some errors are tolerable and where the
data is constantly changing, and can also provide a fast check to
rule out bad inputs.  An additional motivation for studying property
testing is that this area is abundant with fascinating combinatorial
problems. Property testing has recently become an active research
area; a good recent survey is \cite{icm}.

  In this paper, we  address several algorithmic  problems in
the theory of groups and semigroups that may admit sublinear time
solution, at least for ``most" inputs. One of these problems is a
special case of the well-known Whitehead's problem: given two
elements of a free group $F$, find out whether or not one of them
can be taken to the other by an automorphism of $F$. This problem
was solved long time ago by Whitehead himself, but the complexity of
the solution is still a subject of active research.  It is not hard
to show, for example, that those elements (represented by freely
reduced words) which cannot be taken to a free generator (i.e.,
non-primitive elements) can be detected by a sublinear (with respect
to the length of an input element) time algorithm with a negligible
probability of failure; see our Section \ref{back}.

Another problem that we  consider is the word problem. It is fairly
easy to show that testing  sublinear-length subwords of a given
(freely reduced) word $g$ cannot help in deciding whether or not $g
=1$ in $G$ unless $G$ is a free group because one has to at least
test a subword of length about $\frac{1}{2}|g|$. However, with
semigroups the situation is different, so we want to find (natural)
examples of semigroups where the word problem admits a sublinear
time solution for ``most" inputs. One potential source of such
examples is ``positive monoids" associated with groups, i.e.,
monoids generated by group generators, but not their inverses. We
address this problem in Section \ref{wp} for positive monoids
associated with free nilpotent group, with Thompson's group $F$, and
with braid groups. It turns out that of these positive monoids, only
those associated with braid groups admit sublinear-time detecting of
inequality  at least for some pairs of words.

\section{Background: sublinear time algorithms in graph theory}
\label{back}

 There are several interesting contexts where one can test properties in sublinear
time. For example, in \cite{GGR}, the authors focused their
attention on testing various properties of graphs and other
combinatorial objects. In particular, they considered the property
of $k$-colorability. This property is NP-complete to determine
precisely but it is easily testable; more specifically,  one can
distinguish $k$-colorable graphs from those that are $\epsilon$-far
from $k$-colorable in constant time. (Two graphs $G$ and  $H$ on $n$
vertices are $\epsilon$-close if at most $\epsilon n^2$ edges need
to be modified (inserted or deleted) to turn $G$ into $H$.
Otherwise, $G$ and  $H$ are $\epsilon$-far.)

  The work of \cite{GGR} sparked a flurry of other results; in
particular, an interesting line of work was initiated in
\cite{AFKS}, where the authors showed that the property of a graph
being $H$-free (that is, the graph does not contain any copy of $H$
as a subgraph) is easily testable for any constant sized graph $H$.

In general, the area of property testing has been very active, with
a number of property testers suggested for graphs and  other
combinatorial objects, as well as matrices, strings, metric spaces,
etc.

In this paper, we discuss sublinear time property testing in the
context of some particular problems in combinatorial theory of
groups and semigroups. Testing some of these properties amounts to
testing a graph (e.g. the Whitehead graph of a free group element),
and therefore fits in with the original ideas of sublinear time
property testing that come from graph theory. To give an example, we
describe here  a particular property of a free group element that
can be tested in sublinear  time in the length of the input element.

\section{Testing primitivity in a free group}

 Let $F_r$  be a free group of rank $r \ge 2$ with a fixed finite basis $X=\{x_1,\dots,x_r\}$.
An element  $g \in F_r$ is  called {\it primitive} if it is a member
of some free basis of $F_r$.  Or, equivalently, if there is an
automorphism of $F_r$  that takes  $g$  to $x_1$.

A  natural property of a given element  $u \in F_r$ one might want
to test is whether or not $u$ is primitive. We show that for ``most"
inputs, this can be done in time sublinear  in the length of $u$. We
have to note one subtle distinction between what we are going to
show here and what was established in \cite{KSS}. From the results
of \cite{KSS}, it follows that a ``generic" element  $u \in F_r$ is
{\it not} primitive (moreover, its length cannot be decreased by any
automorphism of $F_r$). However, these results are only applicable
if $u$ was chosen uniformly randomly from the set of all (freely
reduced) words of length $\le N$, for some $N$.  Furthermore, given
a particular element  $u \in F_r$, the results of \cite{KSS} do not
allow one to check (in linear time, say) that $u$ is, indeed,
non-primitive.

What we are going to show here is that, after testing a small part
of a cyclically reduced word $u$, one can, for a ``generic" freely
reduced $u$,  tell {\it for sure} (i.e., with a rigorous proof) that
$u$ is not primitive. To explain this, we have to introduce the
Whitehead graph first.

The Whitehead graph $Wh(u)$  of a (cyclically reduced) word $~u \in
F_r$ is obtained as follows. The vertices of this graph correspond
to the elements of the free generating set $X$ and their inverses.
For  each occurrence of a subword  $x_i x_j$ in the  word $u$, there
is an  edge in $Wh(u)$ that connects the vertex $x_i$ to the vertex
$x_j^{-1}$; ~if $u$ has a subword  $x_i x_j^{-1}$, then there is an
edge connecting $x_i$ to $x_j$, etc. ~There is one more edge (the
external edge) included in the definition of the Whitehead  graph:
this is the edge that connects the vertex corresponding to the last
letter of $u$ to the vertex corresponding to the inverse of the
first letter.

 It was observed by Whitehead himself (see also \cite{St99})  that the Whitehead graph
of any cyclically reduced  primitive element of length $>2$ has
either an isolated edge or a cut vertex, i.e., a vertex  that,
having been removed from the graph together with all  incident
edges, increases the number of connected components of the graph.
Obviously, if the Whitehead graph has a Hamiltonian circuit (i.e., a
circuit that contains all vertices of the graph), then it cannot
have a cut vertex. Our test is therefore pretty simple: pick a
random subword $v$ of $u$, of length sublinear in $|u|$, say, of
length  $|u|^\delta$ for some $0 < \delta <1$. It follows  from
results of \cite{KSS} that all possible 2-letter subwords are going
to be present in $v$ with overwhelming probability. Having checked
that (which takes linear time in $|v|$, and therefore sublinear time
in $|u|$), we conclude that the Whitehead graph of $v$ is complete,
hence the Whitehead graph of $u$ has a Hamiltonian circuit, whence
$u$ is not primitive.

We note, in passing, that the problem of detecting a Hamiltonian
circuit in  an  arbitrary given graph is well known to be
computationally hard (in fact, NP-complete) in the worst case
\cite{GJ}, but it is also known to be easy for ``most" graphs (it is
even easy ``on average", see \cite{GS}).

It is an interesting question whether  sublinear time algorithms can
be found for other instances of the Whitehead problem (=automorphic
conjugacy problem) in a free group, so we ask:

\begin{problem}
Let $v \in F_r$ be arbitrary but fixed. Is there a generic subset
$S$ (see our Section \ref{wp}, Definition \ref{defn:as}) of $F_r$
and an algorithm ${\mathcal A}_v$ such that for any $u \in S$, the
algorithm ${\mathcal A}_v$ is able to detect, in time sublinear in
the length of $u$, that  $u$ cannot  be taken to  $v$ by any
automorphism of $F_r$ ?

\end{problem}

\section{The word problem in semigroups}
\label{wp}

If a group (or a semigroup) $G$ is given by a recursive presentation
in terms of generators and defining relators:

$$G = \langle x_1, x_2,\ldots, x_n  \mid r_1, r_2,\dots \rangle,$$

\noindent then the {\it word problem} for $G$ is: given  a word $g =
g(x_1, x_2,\ldots, x_n)$, find out whether or not $g =1$ in $G$. The
word problem is known to have  linear time solution for {\it
hyperbolic groups}.

As we have mentioned in the Introduction, it is fairly easy to show
that testing sublinear-length subwords of a given word $g$ cannot
help in deciding whether or not $g =1$ in a group $G$ unless $G$ is
free and $g$ is freely reduced. Indeed, suppose generators of $G$
satisfy a relation  $r = r(x_1,  \ldots, x_n)=1$. Then, given a
(freely reduced) word $g$ of length $m$, the initial segment of  $g$
of length  $\le \frac{m}{2}$ will have $r$ as a subword with
probability converging to 1 exponentially fast as $m \to \infty$.
Since any cyclic shift of a word representing the identity also
represents the identity, we may assume, without loss of generality,
that our initial segment of  $g$ of length $\le \frac{m}{2}$ ends
with $r$, i.e., it is of the form $ur$. Then, if $g$ is of the form
$urr'u^{-1}$, where $r'$ is any relator in $G$, it represents the
identity. Therefore, examining a subword of length $\le \frac{m}{2}$
of a generic word of length $m$ cannot possibly help to guarantee
that $g \ne 1$ in $G$.

However, with semigroups the situation is different, so we address
here the following, perhaps somewhat vague, problem:

\begin{prob}
Are there natural examples of  semigroups given by generators and
defining relators, where the word problem admits a sublinear time
solution for ``most" inputs?
\end{prob}

Note that the word problem for semigroups has a slightly different
wording (excuse the pun): given  two words $g, h$ in generators of a
semigroup $G$, find out whether or not $g =h$ in $G$.
 Of course, if an algorithm for a sublinear
time solution of the word problem exists, it will only give
``negative"  answers, i. e., $g \ne h$ in $G$. This is similar to
results of \cite{KMSS}, where (generically) linear time solution of
the word problem was offered for several large classes of groups;
their solution, too, gives only ``negative"  answers.

First we have to clarify the meaning of ``most" inputs in this
context. To that  end, we recall the definition of a {\it generic}
set from \cite{KMSS}. The most general and  straightforward
definition is based on the notion of {\it asymptotic density}.

\begin{defn}\label{defn:as}
  Suppose that $T$ is a countable set and that $\ell:T\to \mathbb N$
  is a function (referred to as \emph{length}) such that for every
  $n\in \mathbb N$ the set $\{x\in T: \ell(x)\le n\}$ is finite.  If
  $X\subseteq T$ and $n\ge 0$, we denote $\rho_{\ell}(n,X):=\#\{x\in
  X: \ell(x)\le n\}$ and $\gamma_\ell(n,X)=\#\{x\in X: \ell(x)=n\}$.

  Let $S\subseteq T$. The \emph{asymptotic density} of $S$ in $T$ is
\[
\overline{\rho}_{T,\ell}(S):=\limsup_{n\to\infty}\frac{\#\{x\in S:
  \ell(x)\le n\}}{\#\{x\in T: \ell(x)\le n\}}=\\
\limsup_{n\to\infty}\frac{\rho_\ell(n,S)}{\rho_\ell(n,T)},
\]
where we treat a fraction $\frac{0}{0}$, if it occurs, as $0$.

 If the actual
limit exists, we denote it by $\rho_{T,\ell}(S)$ and call this limit
the \emph{strict asymptotic density} of $S$ in $T$. We say that $S$
is \emph{generic in $T$ with respect to $\ell$} if
$\rho_{T,\ell}(S)=1$.
\end{defn}

 In our situation, $T$ is the set of all words in a given (finite)
alphabet $X = \{x_1, x_2,\ldots, x_n\}$, and $\ell(w), ~w \in T$ is
the usual lexicographic length of     $w$ that we often denote
simply by $|w|$. Thus, given a semigroup $G$ generated by $X$, we
are looking for a generic set $S \subseteq T$ of words such that for
any $g, h \in S$, there is a sublinear time in $n=|g|+|h|$
(probabilistic) algorithm proving that $g \ne h$ in $G$ with
probability $1-\epsilon(n)$, where $\epsilon(n) \to 0$ as $n \to
\infty$.

As we have pointed out in the Introduction, one potential source of
semigroups with the property in question is ``positive monoids"
associated with groups, i.e., monoids generated by group generators,
but not their inverses. For some particular groups, e.g. for braid
groups, Thompson's group,  these monoids have been extensively
studied, and because of very nontrivial combinatorics involved in
these studies, it would be quite interesting to either obtain a
sublinear time algorithm for solving the word problem in these
monoids  or prove that none exists.  Negative results would be
interesting, too, because lower bounds on complexity are always
valuable.

Another important class of positive monoids is associated with free
nilpotent groups; these monoids have a special name of {\it strictly
nilpotent semigroups}, see \cite{semi}. They are called {\it
strictly} nilpotent because there are several other definitions of
nilpotency for semigroups; for a survey on these and on how they are
related to strictly nilpotent semigroups we refer to \cite{semi} or
\cite{Shn}. Here we just say that nilpotent semigroups, under
various definitions, have been extensively studied from many
different perspectives (see e.g. \cite{Olga} or \cite{Shn}).


In the following three subsections, we are going to show that of the
three kinds of positive monoids (associated with free nilpotent
groups, with Thompson's group $F$, and with braid groups), only
those associated with braid groups admit sublinear-time detecting of
inequality at least for some pairs of words.

\subsection{Positive monoid of a free nilpotent group}
\label{nilpotent}

\medskip

Positive monoids of  free nilpotent groups are called  {\it strictly
nilpotent semigroups}, see \cite{semi}. We have to give some
background here because properties of free nilpotent groups are not
as well known these days as properties of braid groups or Thompson's
group are.

Magnus \cite{Magnus} considered the embedding of the free group~$F$
with a free generator set~$X$ into the power series ring with the
same set~$X$ of generators and proved that a~group element which
belongs to $\gamma_c(F)$,  the $c$th term of the lower central
series of the group~$F$, is mapped to a~power series without
non-constant terms of degree less than~$c$. The converse result
(i.e., that any group element which does not belong to
$\gamma_c(F)$, is mapped to a~power series with some non-constant
terms of degree less than~$c$) appeared to be quite difficult to
prove. Probably the first full and correct proof was given by Chen,
Fox and Lyndon in~\cite{CFL}. They considered the free group ring
instead of the power series ring and proved that
$$
\gamma_c(F)=(\Delta_F^c+1)\cap F,
$$

\noindent where  $\Delta_F$ is the {\it augmentation ideal} of the
free group ring $ZF$, i.e., the kernel of the natural
``augmentation" homomorphism $\varepsilon_F\colon ZF\to Z$ that
takes all elements of $F$ to 1.

Now let $M_c$ denote the positive monoid of the free nilpotent group
$F/\gamma_{c+1}(F)$ of class $c$, where $F$ is a free group of rank
$r \ge 2$ with a free generator set~$X$. We do not include the rank
$r$ in the notation because our results in this section are
independent of $r$. We have:

\begin{lemma}\label{semi}
Elements of $M_c$ satisfy all  identities of the form $a_c=b_c$,
where $a_c$ and $b_c$ are words in~$X$ such that
$(a_c-b_c)\in\Delta^{c+1}$.

\end{lemma}

To prove the main result of this section, we will need to combine
this lemma with the following result due to A.~I.~Mal'cev~\cite{15}
and, independently, to B.~Neumann and T.~Taylor \cite{16}. To better
tailor (no pun intended) this result to our needs, we give it here
in a weaker form.

\begin{lemma}\label{Malcev}
Let
$$ u_0=x,\ \ v_0=y;\quad u_{n+1}=u_nv_n,\ \
v_{n+1}=v_nu_n,
$$
where $x$, $y$ are arbitrary elements of a free group $F$. Then
$u_c=v_c $ modulo $\gamma_{c+1}(F)$.
\end{lemma}

By combining Lemma \ref{semi} and  Lemma \ref{Malcev}, we get

\begin{proposition}
For any two  positive words   $w_1$  and $w_2$ of  length $n$ in an
alphabet $X$, there are positive words   $z_1$  and $z_2$ of lengths
 $\le (n-1) \cdot 2^c$ such that $w_1 z_1 = w_2 z_2$ in $M_c$.
\end{proposition}

In particular, given two positive words of length $L$ one cannot
tell that they are not equal in $M_c$ by just inspecting the
prefixes of length $\le \frac{L}{2^c}$, i.e., there is at least no
obvious sublinear time algorithm for detecting inequality in $M_c$.

\begin{proof}
Construct the Mal'cev-Neumann-Taylor sequence of words starting with
$u_0=w_1, ~v_0=w_2$. Then $u_c$ has $w_1$ as a prefix, $v_c$ has
$w_2$ as a prefix, $u_c=v_c$ in $M_c$, and the length of both $u_c$
and $v_c$ is  ~$n \cdot 2^c$.
\end{proof}

\subsection{Positive monoid of Thompson's group $F$}
\label{Thompson}

Thompson's group $F$ is well known in many areas of mathematics,
including algebra, geometry, and analysis.  For a survey on various
properties of Thompson's group, we refer to \cite{CFP}. This group
has the following nice presentation in terms of generators and
defining relations:

$$F = \langle x_0 , x_1 , x_2 , \ldots \mid  x_k x_i =
x_i x_{k+1} ~ (k>i) \rangle.$$

Since all defining relators in this presentation are pairs of
positive words, we can consider the positive  monoid associated with
this presentation; denote it by $F^+$.

We note that the above (infinite) presentation  allows for a
convenient normal form. We do not really need it in this paper, but
we describe it  here anyway. The classical normal form of an element
of Thompson's group is a word of the form
$$x_{i_1} \ldots x_{i_s} x_{j_t}^{-1} \ldots x_{j_1}^{-1},$$
\noindent such that the following two conditions are satisfied:
\begin{enumerate}
 \item[{\bf (NF1)}] $i_1 \le ...\le i_s$ and $j_1\le \ldots \le j_t$
 \item[{\bf (NF2)}] if both $x_i$ and $x_i^{-1}$ occur, then either
$x_{i+1}$ or $x_{i+1}^{-1}$  occurs, too.
\end{enumerate}

Now we get to the point of this section.

\begin{proposition} \label{guba}
For any two  positive words   $w_1$  and $w_2$ of  lengths $m$ and
$n$, respectively, in the alphabet $X=\{x_0 , x_1 , x_2 , \ldots
\}$, there are positive words $z_1$  and $z_2$ of lengths $n$ and
$m$, respectively, such that $w_1 z_1 = w_2 z_2$ in Thompson's group
$F$.
\end{proposition}

The following elegant and simple proof is due to Victor Guba.

\begin{proof}
Construct the following van Kampen diagram (see e.g. \cite{LS} for
the definition of a van Kampen diagram). On a square lattice, mark
one point as the origin. Starting at the origin and going to the
right, write the word   $w_1$ by marking edges of the lattice by the
letters of $w_1$, read left to right. Then, starting at the origin
and going up, write the word   $w_2$ by marking edges of the lattice
by the letters of $w_2$, read left to right.

Now start marking edges of the lattice inside the rectangle built on
segments of length $m$ (horizontally) and $n$ (vertically)
corresponding to the words   $w_1$  and $w_2$, as follows. All
horizontal edges in the lattice are directed from left to right, and
all vertical edges are directed from bottom to top. Then, suppose a
single square cell of the lattice has:

\begin{itemize}

\item $x_i$ on the lower edge and $x_i$ on the left edge. Then we
mark the upper edge and the right edge of this cell with the same
$x_i$. This cell now corresponds to the relation $x_ix_i=x_ix_i$.

\item $x_i$ on the lower edge and $x_j$ on the left edge, where
$i<j$. Then we mark the upper edge of this cell with $x_i$, and the
right edge with $x_{j+1}$. This cell now corresponds to the relation
$x_jx_ix_{j+1}^{-1}x_i^{-1}=1$,  or  $x_jx_i=x_ix_{j+1}$.

\item $x_i$ on the lower edge and $x_j$ on the left edge, where
$i>j$. Then we mark the upper edge of this cell with $x_{i+1}$, and
the right edge with $x_{j}$. This cell now corresponds to the
relation $x_jx_{i+1}x_{j}^{-1}x_i^{-1}=1$,  or
$x_jx_{i+1}=x_ix_{j}$.

\end{itemize}

After all edges of the rectangle built on segments corresponding to
the words   $w_1$  and $w_2$ are marked, we read a relation of the
form $w_2 u_1 u_2^{-1} w_1^{-1}=1$,  or $w_2 u_1=w_1u_2$,  off the
edges of this rectangle. Here the length of $u_1$ is $m$ and the
length of $u_2$ is $n$. This completes the proof.

\end{proof}

\begin{example} If $w_1=x_1x_2$  and $w_2=x_3x_5$, this method
gives $w_1 x_5x_7 = w_2 x_1x_2$.
\end{example}

Proposition \ref{guba} implies, in particular, that it is impossible
to tell that two positive words of length  $L$ in the alphabet
$X=\{x_0 , x_1 , x_2 , \ldots \}$ are not equal in Thompson's group
$F$ by inspecting their initial segments of length $\le
\frac{L}{2}$, i.e., there is at least no such straightforward
sublinear time algorithm for detecting inequality in $F^+$.

\subsection{Positive braid monoids}
\label{braids}

Braid groups need no introduction; we just refer to the monograph
\cite{B} for  background. Some notation has to be  recalled though.
We denote the braid group on $n$ strands by $B_n$; this group has a
standard presentation
$$\langle \sigma_1, ..., \sigma_{n-1}| ~
\sigma\sb i \sigma\sb j = \sigma\sb j \sigma\sb i ~{\rm if}
 ~\vert i-j\vert  > 1;
~\sigma\sb i \sigma\sb {i+1} \sigma\sb i = \sigma\sb {i+1} \sigma\sb
i \sigma\sb {i+1}
 ~{\rm for}\ 1 \le i \le n-2 \rangle.$$
We shall call elements of $B_n$ {\it braids}, as opposed to {\it
braid words} that are elements of the ambient free group on
$\sigma_1, ..., \sigma_{n-1}$.

Since all defining relators of a braid group are positive words, we
can consider the positive braid monoid; denote it by $B_n^+$.

It turns out that, in contrast to the situation with positive
monoids  $M_c$ and  $F^+$ considered in two previous sections of
this paper, for at least some pairs of positive words in $B_n^+$
there is a sublinear time test for inequality. The following
proposition follows from the results of \cite{AutordDehornoy}; in
particular, from the proof of their Proposition 2.9.

\begin{proposition} \label{braids}
Let $w_1 = \sigma_1 \sigma_3  \cdots \sigma_{2m-1}$, ~~$w_2 =
\sigma_{2m} \sigma_{2m-2} \cdots \sigma_{2}$. Suppose $w_1 u = w_2
v$ for some $u, v \in B_n^+, ~n \ge 2m$. Then $|u|, |v| = 2m^2$.
\end{proposition}

Thus, in particular, if one has two positive braid words of length
$L$, where one of them starts with  $\sigma_1 \sigma_3  \cdots
\sigma_{2k-1}$, the other one starts with  $\sigma_{2k}
\sigma_{2k-2} \cdots \sigma_{2}$,  and $k \ge \sqrt{L}$, then these
braid words are not equal in $B_n^+, ~n \ge 2k$.

Of course, this is just a very special example where a sublinear
time algorithm can detect inequality of two words in $B_n^+$, so the
interesting question is whether examples of this sort are
``generic". We therefore ask:

\begin{problem}
Is there a generic subset $S$ (in the sense of Definition
\ref{defn:as}) of $B_n^+$ and a number $\epsilon > 0$    such that
for any two words $w_1, w_2$ of length $k$ representing elements of
$S$, the minimum length of words $u, v$ such that $w_1u =  w_2v$, is
greater than $k^{(1+\epsilon)}$ ?\\

\end{problem}

\noindent {\it Acknowledgement.} The author is grateful to Victor
Guba and Patrick Dehornoy for helpful discussions.

{}

\end{document}